\documentclass[a4paper]{article}
\usepackage{amsmath}
\usepackage{amsfonts}
\usepackage{amsthm}
\usepackage{amssymb}
\usepackage[USenglish]{babel}
\usepackage[utf8]{inputenc} 
\usepackage[T1]{fontenc}
\usepackage{times}
\usepackage{cite}

\theoremstyle{remark}
\newtheorem{remark}{Remark}[section]

\begin{document}

\title{The $\operatorname{SW}(3/2,2)$ superconformal algebra via a Quantum Hamiltonian Reduction of $\operatorname{osp}(3|2)$}

\author{L\'azaro O. Rodr\'iguez D\'iaz}

\date{}

\maketitle

\begin{abstract}
We prove that the family of non-linear $W$-algebras $\operatorname{SW}(3/2,2)$ which are extensions of the $N=1$ superconformal algebra by a primary supercurrent of conformal weight $2$ can be realized as a quantum Hamiltonian reduction of the Lie superalgebra $\operatorname{osp}(3|2)$. In consequence we obtain an explicit free field realization of the algebra in terms of the screening operators.  At central charge $c=12$ the $\operatorname{SW}(3/2,2)$ superconformal algebra corresponds to the superconformal algebra associated to sigma models based on eight-dimensional manifolds with special holonomy $\operatorname{Spin(7)}$, i.e., the Shatashvili-Vafa $\operatorname{Spin(7)}$ superconformal algebra. 
\end{abstract}

\section{Introduction}

The $\operatorname{SW}(3/2,2)$ superconformal algebra was first constructed in \cite{FarrilSchrans91, FarrilSchrans92}. It is a one parameter family of non-linear $W$ algebras, parametrized by the central charge $c$. This algebra is an extension of the $N=1$ superconformal algebra $\{G,L\}$ by a primary supercurrent of conformal weight $2$, that is, by two primary fields $W$ and $U$ of conformal weight $2$ and $5/2$ respectively, that satisfy:
\begin{align*}
[G_{\lambda}W]=U,\qquad [G_{\lambda}U]=(\partial+4\lambda)W.\nonumber
\end{align*}
The other $\lambda$-brackets are:
\begin{align}\label{eq:bracket of W with W}
[W_{\lambda} W] = &\tfrac{c}{12}\lambda^{3}
+ \left(
  2L +{\tfrac{2(6 + 5c)}
   {{\sqrt{15 - c}}{\sqrt{21 + 4c}}}} W \right)\lambda + \nonumber\\
& +\partial L +{\tfrac{6 + 5c}
   {{\sqrt{15 - c}}\,{\sqrt{21 + 4c}}}} \partial W,\nonumber
\end{align}
\begin{align}  
[W_\lambda U]=&(-\tfrac{3}{2})G\lambda^{2}+\left(\tfrac{6+ 5c}{\sqrt{15-c}\sqrt{21+4c}}U-\partial G\right)\lambda \nonumber \\
&+\tfrac{(15-c)}{(21+4c)}\partial^2 G +\tfrac{(-2)\sqrt{15-c}}{\sqrt{21+4c}}\partial U +\tfrac{(-54)}{(21+4c)}:LG:\nonumber\\
&+\tfrac{(-54)}{\sqrt{15-c}\sqrt{21+4c}}:WG:,\nonumber
\end{align}
\begin{align}
 [U_\lambda U]=&-\tfrac{c}{12}\lambda^{4}+\left(-5L-\tfrac{2(6+5c)}{\sqrt{15-c}\sqrt{21+4c}} W\right)\lambda^{2}\nonumber \\
&+\left(-5\partial L-\tfrac{2(6+5c)}{\sqrt{15-c}\sqrt{21+4c}}\partial W\right)\lambda \nonumber \\
&-\tfrac{6(c+3)}{(21+4c)}\partial^{2}L+\tfrac{3(6-c)}{\sqrt{15-c}\sqrt{21+4c}}\partial^{2}W\nonumber\\
&-\tfrac{108}{(21+4c)}:LL:-\tfrac{108}{\sqrt{15-c}\sqrt{21+4c}}:LW:\nonumber\\
&+\tfrac{27}{(21+4c)}:G\partial G: + \tfrac{54}{\sqrt{15-c}\sqrt{21+4c}}:GU:\nonumber
\end{align}
The last two $\lambda$-brackets contain product of fields, therefore the non-linearity of the algebra.  At central charge $c=12$ the $\operatorname{SW}(3/2,2)$ superconformal algebra corresponds to the superconformal algebra associated to sigma models based on eight-dimensional manifolds with special holonomy $\operatorname{Spin(7)}$, i.e., the Shatashvili-Vafa $\operatorname{Spin(7)}$ superconformal algebra \cite{Shatashvili-Vafa95}; the generator $W$ comes from the invariant four form on a $\operatorname{Spin(7)}$ manifold. The $\operatorname{SW}(3/2,2)$ algebra appeared in \cite{Mallwitz95} as the algebra of conserved currents of the quantized supersymmetric Toda theory corresponding to $\operatorname{osp}(3|2)$. A classical version of this result was obtained in \cite{NoharaMohri91}.

In section \ref{sec-reduction} we obtain the $\operatorname{SW}(3/2,2)$ superconformal algebra via a quantum Hamiltonian reduction of the Lie superalgebra $\operatorname{osp}(3|2)$, as a corollary we get a free field realization of the algebra in terms of two free bosons and two free fermions, as well as an explicit description of the screening operators.

We perform the quantum Hamiltonian reduction using the general framework developed in\cite{KacWakimoto04}. Let $\mathfrak{g}$ be a simple Lie superalgebra with a non-degenerate even supersymmetric invariant bilinear form $(.|.)$,\; $k\in\mathbb{C}$ and $x, f \in \mathfrak{g}$ a pair of even elements such that $\operatorname{ad}x$ is diagonalizable on $\mathfrak{g}$ with half-integer eigenvalues $\left(\mathfrak{g}=\oplus_{j\in\frac{1}{2}\mathbb{Z}}\mathfrak{g}_{j}\right)$, $f\in \mathfrak{g}_{-1}$ a nilpotent element and the eigenvalues of $\operatorname{ad}x$ on the centralizer $\mathfrak{g}^{f}$ of $f$ in $\mathfrak{g}$ are non-positive. To this data is associated a homology complex 
\begin{align*}
C(\mathfrak{g},x,f,k):=(V_{k}(\mathfrak{g})\otimes F(\mathfrak{g}_{+}\oplus\mathfrak{g}_{+}^{*})\otimes F(\mathfrak{g}_{1/2}), d_{0})
\end{align*}
 where $V_{k}(\mathfrak{g})$ is the universal affine vertex superalgebra of level $k$ associated to $\mathfrak{g}$, $F(\mathfrak{g}_{+}\oplus\mathfrak{g}_{+}^{*})$ is the vertex superalgebra of charged free superfermions associated to $\mathfrak{g}_{+}\oplus\mathfrak{g}_{+}^{*}$ with the reversed parities, $F(\mathfrak{g}_{1/2})$ is the neutral free superfermion vertex superalgebra associated to $\mathfrak{g}_{1/2}$ and $d_{0}$ a is suitable differential. The $W$-algebra $W_{k}(\mathfrak{g},x,f)$ is defined as the homology  $W_{k}(\mathfrak{g},x,f):=H_{\bullet}(C(\mathfrak{g},x,f,k))$ of the complex together with the vertex superalgebra structure induced from $C(\mathfrak{g},x,f,k)$. In fact, the structure theorem \cite[Theorem 4.1]{KacWakimoto04} proved $W_{k}(\mathfrak{g},x,f)$ is the $0$-th homology of the complex and that it is strongly generated by fields $J^{\{u_{i}\}}$, where $\{u_{1}, \cdots, u_{n}\}$ is a basis of $\mathfrak{g}^{f}$, moreover, as the complex has the formality property it is possible to compute in principle the generators $J^{\{u_{i}\}}$. For more details and the notation used in the present paper we refer the reader to \cite{KacWakimoto04} or \cite[Section 2]{HeluaniRodriguez15}.

Observe that a quantum Hamiltonian reduction of $\operatorname{osp}(3|2)$ with respect to the minimal nilpotent element $f$, i.e., $f$ is a root vector associated to an even highest root of $\mathfrak{g}$ produces the $N=3$ superconformal algebra after being tensored with one fermion \cite{KacWakimoto04}.

In \cite{GepnerNoyvert01} the study of the representation theory of $\operatorname{SW}(3/2,2)$ was addressed, however we still lack a character formulae. Combining recent advances in the representation theory of affine Lie superalgebras \cite{GorelikKac15}, in particular the Kac-Wakimoto character formula for $\widehat{\operatorname{osp}}(3|2)$, we can obtain a character formula for $\operatorname{SW}(3/2,2)$ once we apply the Kac-Frenkel-Wakimoto functor to our construction; these will appear in a forthcoming paper. Recently \cite{Kachru16} was observed that if we decompose the elliptic genus of a $\operatorname{Spin(7)}$ holonomy manifold in terms of conjectural characters of the Shatashvili-Vafa $\operatorname{Spin(7)}$ superconformal algebra then we obtain a similar phenomena as observed in Mathieu moonshine.

\section{Quantum Hamiltonian reduction of $\operatorname{osp}(3|2)$}\label{sec-reduction}

In this section we prove that the $SW(\frac{3}{2},2)$ superconformal $W$-algebra can be realized as a quantum Hamiltonian reduction of $\operatorname{osp}(3|2)$. As a direct consequence we obtain a free field realization in terms of two free bosons and two free fermions along with a description of the screening operators. 

The Lie superalgebra $\operatorname{osp}(3|2)$  is a simple Lie superalgebra of rank $2$ and dimension $12$, see \cite{Kac77}. It has non-vanishing Killing form and dual Coxeter number $h^{\vee}=\frac{1}{2}$. As a contragradient Lie superalgebra $\mathfrak{g}=\operatorname{osp}(3|2)$ is associated to the Cartan matrix $A=(a_{ij})_{i,j}$ and $\tau=\{1\}$
\begin{eqnarray*}\label{matrixCartan}
{(a_{ij})}_{i,j=1}^{2}= \left( \begin{array}{cc}
0 & \frac{1}{2}  \\
\frac{1}{2} & -\frac{1}{2} \\
\end{array} \right).
\end{eqnarray*}
It has generators $\{h_{1},h_{2}, e_{1},e_{2},f_{1},f_{2}\}$, all being even except $e_1,f_1$. Let's define
\begin{gather}
[e_1,e_2]=:e_{12}, \quad [e_2,e_{12}]=:e_{122},\quad [e_1,e_{122}]=:e_{1122},\nonumber\\
[f_1,f_2]=:f_{12},\quad [f_2,f_{12}]=:f_{122},\quad [f_1,f_{122}]=:f_{1122},\nonumber
\end{gather}
and denote by $\Pi=\{\alpha_{1},\alpha_{2}\}$ its simple roots.\\
\noindent
To carry out the quantum Hamiltonian reduction we use the pair $(x,f)$:
\begin{align}
x:=h_{1}-h_{2},\quad f:=f_{2}+f_{1122}.\nonumber
\end{align}
and by non-degenerate even supersymmetric invariant bilinear form $(.|.)$ we take the Killing form:
\begin{gather}
(h_i,h_j)=a_{ij},\quad(e_i,f_j)=\delta_{ij},\quad (e_{12},f_{12})=-(f_{12},e_{12})=\tfrac{1}{2},\nonumber\\
(e_{122},f_{122})=-(f_{122},e_{122})=\tfrac{1}{4},\quad (e_{1122},f_{1122})=(f_{1122},e_{1122})=-\tfrac{1}{4}.\nonumber
\end{gather}
The eigenspace decomposition $\mathfrak{g}=\oplus_{j\in\frac{1}{2}\mathbb{Z}}\mathfrak{g}_{j}$ of the algebra with respect to ad\,$x$ is as follows:
\[  \begin{array}{ccccccc}
\mathfrak{g}_{-3/2} & \mathfrak{g}_{-1} & \mathfrak{g}_{-1/2} & \mathfrak{g}_{0} & \mathfrak{g}_{1/2} & \mathfrak{g}_{1} & \mathfrak{g}_{3/2} \\
f_{122} & f_{2} & e_{1} & h_{1} & f_{1} & e_{2} & e_{122} \\
 & f_{1122} & f_{12} & h_{2} & e_{12} & e_{1122} & \\
 \end{array} \]
where $\mathfrak{g}_{j}:=\{a\in\mathfrak{g}\;|\; [x,a]=j\,a\}$. Denote by $\mathfrak{g}^{f}$ the centralizer of $f$ in $\mathfrak{g}$, then we have $\mathfrak{g}^{f}=\mathfrak{g}_{-1/2}^{f}\oplus \mathfrak{g}_{-1}^{f}\oplus \mathfrak{g}_{-3/2}^{f}$ with $\operatorname{dim}\mathfrak{g}_{-1/2}^{f}=1$ , $\mathfrak{g}_{-1}^{f}\simeq \mathfrak{g}_{-1}$ and $\mathfrak{g}_{-3/2}^{f}\simeq \mathfrak{g}_{-3/2}$.\\
\noindent 
Denote by $F(\mathfrak{g}_{1/2})$ the neutral free superfermion vertex superalgebra associated to the vector superspace $\mathfrak{g}_{1/2}$
with the even skew-supersymmetric non-degenerate bilinear form:
\begin{align}
\left<a|b\right>_{ne}=(f|[a,b]), \quad a, b \in \mathfrak{g}_{1/2}.\nonumber
\end{align} 
Let's denote by $\Phi_{-1}:=f_{1}$ and $\Phi_{12}:=e_{12}$  the free neutral fermions. They satisfy the following non-zero $\lambda$-brackets:
\begin{align*}
[{\Phi_{-1}}_{\lambda}\Phi_{12}]=\left<\Phi_{-1}|\Phi_{12}\right>_{ne}=\frac{1}{2},\\
[{\Phi_{12}}_{\lambda}\Phi_{12}]=\left<\Phi_{12}|\Phi_{12}\right>_{ne}=-\frac{1}{4},
\end{align*} 
and the dual free neutral fermions with respect to $\left<.|.\right>_{ne}$ are:
\begin{gather*}
\Phi^{-1}=\Phi_{-1}+2\Phi_{12},\quad \Phi^{12}=2\Phi_{-1}.
\end{gather*}
Let $V_{k+\frac{1}{2}}(\mathfrak{g})$ denote the affine vertex algebra of level $k+\tfrac{1}{2}$ associated to $\mathfrak{g}$, and denote its currents by $J^{(v)}$.
\begin{remark}\label{rmk-free field}
 Performing the quantum Hamiltonian reduction we obtain $W_{k}(\mathfrak{g},x,f)$ as a subalgebra of $V_{k+\frac{1}{2}}(\mathfrak{g}_{\le})\otimes F(\mathfrak{g}_{1/2})$, where $\mathfrak{g}_{\le}:=\bigoplus_{j\le 0}\mathfrak{g}_{j}$. Moreover, as $\mathfrak{g}_{0}$ equals the Cartan subalgebra $\mathfrak{h}$ of $\mathfrak{g}$ the canonical homomorphism $\mathfrak{g}_{\le}\rightarrow\mathfrak{g}_{0}$ induces a free field realization of $W_{k}(\mathfrak{g},x,f)$ inside $V_{k+\frac{1}{2}}(\mathfrak{h})\otimes F(\mathfrak{g}_{1/2})$, see \cite[Remark 2.3]{HeluaniRodriguez15}.
\end{remark}
Fix the basis $\{h_{1},h_{2},e_{1},f_{12},f_{2},f_{1122},f_{122}\}$  of $\mathfrak{g}_{\le}$ compatible with the $\frac{1}{2}\mathbb{Z}$ and $\mathbb{Z}_{2}$ gradation of $\mathfrak{g}$. Fix also the basis $\{e_{1}-f_{12}, f_{2}, f_{1122}, f_{122}\}$ of $\mathfrak{g}^{f}$ compatible with its $\frac{1}{2}\mathbb{Z}$ gradation and let's denote by $J^{\{e_{1}-f_{12}\}}, J^{\{f_{2}\}}, J^{\{f_{1122}\}}, J^{\{f_{122}\}}$ the corresponding fields of conformal weights $\tfrac{3}{2}, 2, 2, \tfrac{5}{2}$ given by Kac-Wakimoto \cite[Theorem 4.1]{KacWakimoto04} that strongly generate $W_{k}(\mathfrak{g},x,f)$.\\
\noindent 
The generating field in conformal weight $\tfrac{3}{2}$ can be computed using a formula obtained in \cite[Theorem 2.1]{KacWakimoto04}:   
\begin{align*}
J^{\{e_{1}-f_{12}\}}=&J^{(e_{1})}-J^{(f_{12})}+:\Phi^{-1}J^{(h_1)}:-\tfrac{1}{2}:\Phi^{12}J^{(h_1)}:-\tfrac{1}{2}:\Phi^{12}J^{(h_2)}:\\
&-(\tfrac{1}{2}+k)\partial\Phi^{-1}-\tfrac{1}{2}k\partial\Phi^{12}.
\end{align*}
With regard to the other fields we need to work a bit more and exploit the formality of the complex to get an explicit expression. However, as it is possible to recover all the fields in the $SW(\frac{3}{2},2)$ algebra from the generators in conformal weight $\frac{3}{2}$ and $2$, we only need to compute $J^{\{f_{2}\}}$ and $J^{\{f_{1122}\}}$:
\begin{align*}
J^{\{f_2\}}=&J^{(f_2)}-\tfrac{1}{2}:\Phi^{12}J^{(e_1)}:+:\Phi^{-1}J^{(f_{12})}:+:J^{(h_2)}J^{(h_2)}:\\
&+\tfrac{1}{2}:\Phi^{-1}\Phi^{12}J^{(h_1)}:+\tfrac{1}{2}:\Phi^{-1}\Phi^{12}J^{(h_2)}:+\tfrac{1}{2}\left(1+4k\right)\partial J^{(h_2)}\\
&+\tfrac{1}{2}k:\Phi^{-1}\partial\Phi^{12}:,
\end{align*}
\begin{align*}
J^{\{f_{1122}\}}=&J^{(f_{1122})}+\tfrac{1}{2}:J^{(f_{12})}\Phi^{12}:+(-2):J^{(h_1)}J^{(h_2)}:+(-1):J^{(h_1)}J^{(h_1)}:\\
&+(-1):J^{(h_2)}J^{(h_2)}:+(-\tfrac{1}{2}):\Phi^{-1}\Phi^{12}J^{(h_1)}:+(-\tfrac{1}{2}):\Phi^{-1}\Phi^{12}J^{(h_2)}:\\
&+(-k)\partial J^{(h_1)}+(-k)\partial J^{(h_2)}+\tfrac{1}{8}\left(-1-2k\right) :\partial\Phi^{-1}\Phi^{12}:\\
&-\tfrac{1}{8}\left(-1+2k\right):\Phi^{-1}\partial\Phi^{12}:+\tfrac{1}{16}\left(1+2k\right):\partial\Phi^{12}\Phi^{12}:.
\end{align*}
The central charge of the Virasoro field of $W_{k}(\mathfrak{g},x,f)$ as a function of the level $k$ is given by formula $c(k)=6+18k$, see \cite{KacWakimoto04}.\\
\noindent
We proceed now to get the basis $\{L, G, W, U\}$ of the algebra as presented in the introduction. We are looking for a field $G$ such that $\{G, L:=\tfrac{1}{2}{G}_{(0)}G\}$ generate an $N=1$ superconformal algebra of central charge $c(k)$, then:
\begin{align*}
G:=a\left(J^{\{e_{1}-f_{12}\}}\right)
\end{align*}
\noindent
where $a=\frac{2}{\sqrt{\left(-1-2k\right)}}$.

The field $W=a_{1}J^{\{f_2\}}+a_{2}J^{\{f_{1122}\}}$ of conformal weight $2$, get fixed by the conditions:
\begin{gather*}
{G}_{(j)}W=0, \;\;j>0, \quad {W}_{(3)}W=\frac{c(k)}{12},
\end{gather*}
\noindent
that is,
\begin{align*}
a_{1}&=\tfrac{2\sqrt{1-2k}\sqrt{5+8k}}{5+18k+16k^{2}},\\
a_{2}&=\left(\tfrac{2+8k}{2k-1}\right)\tfrac{2\sqrt{1-2k}\sqrt{5+8k}}{5+18k+16k^{2}}.
\end{align*}
Then the explicit formulas of all the generators of $W_{k}(\mathfrak{g},x,f)$ as a subalgebra of $V_{k+\frac{1}{2}}(\mathfrak{g}_{\le})\otimes F(\mathfrak{g}_{1/2})$ look as follows:
\begin{align*}
G=&\tfrac{2}{\sqrt{\left(-1-2k\right)}}\left(J^{(e_{1})}-J^{(f_{12)}}+:\Phi^{-1}J^{(h_1)}:-\tfrac{1}{2}:\Phi^{12}J^{(h_1)}:\right.\\
&\left.-\tfrac{1}{2}:\Phi^{12}J^{(h_2)}:-(\tfrac{1}{2}+k)\partial\Phi^{-1}-\tfrac{1}{2}k\partial\Phi^{12}\right),
\end{align*}
\begin{align*}
L=&\left(-\tfrac{2}{1+2k}\right)\left(J^{(f_2)}+J^{(f_{1122})}-:J^{(h_1)}J^{(h_1)}:
-2:J^{(h_1)}J^{(h_2)}:\right.\\
& \left.+:\Phi^{-1}J^{(f_{12})}:+\left(\tfrac{1+2k}{8}\right):\Phi^{-1}\partial \Phi^{12}:+(-\frac{1}{2}):\Phi^{12}J^{(e_1)}:\right. \\
& \left.+(-\tfrac{1}{2}):\Phi^{12}J^{(f_{12})}:+\left(-\tfrac{1+2k}{8}\right):\partial \Phi^{-1}\Phi^{12}:+ \left(\tfrac{1+2k}{16}\right):\partial \Phi^{12}\Phi^{12}:\right.\\
&\left.+(-1)\partial J^{(h_1)}+\left(\tfrac{1+2k}{2}\right)\partial J^{(h_2)}\right),
\end{align*}
\begin{align*}
W=&\tfrac{2\sqrt{1-2k}\sqrt{5+8k}}{5+18k+16k^{2}}\left(
J^{(f_{2})}+\tfrac{(2+8 k)}{-1+2 k}J^{(f_{1122})}+\tfrac{(2+8 k)}{1-2 k}:J^{(h_{1})}J^{(h_{1})}:\right.\\
&\left.+\tfrac{(4+16 k)}{1-2 k}:J^{(h_{1})}J^{(h_{2})}:+\tfrac{(3+6 k)}{1-2 k}:J^{(h_{2})}J^{(h_{2})}:+:\Phi^{-1}J^{f_{12}}:\right.\\
&\left.+\tfrac{(3+6 k)}{2-4 k}:\Phi^{-1}\Phi^{12}J^{(h_{1})}:+\tfrac{(3+6 k)}{2-4 k}:\Phi^{-1}\Phi^{12}J^{(h_{2})}:\right.\\
&\left.+\tfrac{(-1-2 k)}{4}:\Phi^{-1}\partial\Phi^{12}:-\tfrac{1}{2}:\Phi^{12}J^{(e_{1})}:+\tfrac{(1+4 k)}{1-2 k}:\Phi^{12}J^{(f_{12})}:\right.\\
&\left.+\tfrac{\left(1+6 k+8 k^2\right)}{4-8 k}:\partial\Phi^{-1}\Phi^{12}:+\tfrac{\left(1+6 k+8 k^2\right)}{-8+16 k}:\partial\Phi^{12}\Phi^{12}:\right.\\
&\left.-\tfrac{2 k (1+4 k)}{-1+2 k}\partial J^{(h_{1})}+\tfrac{\left(1+6 k+8 k^2\right)}{2-4 k}\partial J^{(h_2)}\right),
\end{align*}
\begin{align*}
U=&\left(\tfrac{12 i \sqrt{1+3 k} }{\sqrt{1-2 k} \sqrt{5+8 k} \sqrt{1+5 k+6 k^2}}\right)\left(J^{(f_{122})}-:J^{(h_{1})}J^{(e_{1})}:-:J^{(h_{2})}J^{(e_{1})}:\right.\\
&\left.-:J^{(h_{2})}J^{(f_{12})}:+:\Phi^{-1}J^{(f_{1122})}:-:\Phi^{-1}J^{(h_{1})}J^{(h_{1})}:\right.\\
&\left.-:\Phi^{-1}J^{(h_{1})}J^{(h_{2})}:-\tfrac{1}{4}:\Phi^{-1}\Phi^{12}J^{(e_{1})}:-\tfrac{1}{2}:\Phi^{-1}\Phi^{12}J^{(f_{12})}:\right.\\
&\left.+\tfrac{(1+2 k)}{16}:\Phi^{-1}\partial\Phi^{12}\Phi^{12}:+\tfrac{(1-2 k)}{6}:\Phi^{-1}\partial J^{(h_{1})}:+\tfrac{1}{4}:\Phi^{12}J^{(f_{2})}:\right.\\
&\left.-\tfrac{1}{2}:\Phi^{12}J^{(h_{1})}J^{(h_{2})}:-\tfrac{1}{2} :\Phi^{12}J^{(h_{2})}J^{(h_{2})}:+\tfrac{(-1-4 k)}{12}:\Phi^{12}\partial J^{(h_{1})}:\right.\\
&\left.+\tfrac{(-1-4 k)}{12}:\Phi^{12}\partial J^{(h_{2})}:+\tfrac{(1+4 k)}{6}:\partial\Phi^{-1}J^{(h_{1})}:+\left(\tfrac{1}{2}+k\right):\partial\Phi^{-1}J^{(h_{2})}:\right.\\
&\left.+\tfrac{(1+2 k)}{8}:\partial\Phi^{-1}\Phi^{-1}\Phi^{12}:+\tfrac{(-1-4 k)}{12}:\partial\Phi^{12}J^{(h_{1})}:\right.\\
&\left.+\left(-\tfrac{1}{12}-\tfrac{5 k}{6}\right) :\partial\Phi^{12}J^{(h_{2})}:+\tfrac{(1-2 k)}{6}\partial J^{(e_{1})}:+\tfrac{(-1-4 k)}{6}\partial J^{(f_{12})}:\right.\\
&\left.+\tfrac{\left(1+6 k+8 k^2\right)}{24}\partial^{2}\Phi^{-1}-\tfrac{k (1+4 k)}{12}\partial^{2}\Phi^{12}\right).
\end{align*}
\noindent
Using Thielemans's software \cite{Thielemans91} we have checked that the $\lambda$-brackets of the fields above coincide with the $\lambda$-brackets of the $SW(\frac{3}{2},2)$ algebra with central charge $c(k)=6+18k$. The Shatashvili-Vafa $\operatorname{Spin}(7)$ superconformal algebra corresponds to $c=12$, that is, $k=1/3$.\\
\newpage
\noindent
\emph{Free field realization}
\noindent

In order to get the free field realization of $W_{k}(\mathfrak{g},x,f)$ inside $V_{k+\frac{1}{2}}(\mathfrak{h})\otimes F(\mathfrak{g}_{1/2})$  we only need to eliminate the terms that contain a $J^{(v)}$ with $v\in\mathfrak{g}_{\le}\backslash \mathfrak{g}_{0}$, that is, the terms containing $J^{(e)}$'s and $J^{(f)}$'s (see Remark \ref{rmk-free field}):
\begin{align*}
G=&\tfrac{2}{\sqrt{\left(-1-2k\right)}}\left(:\Phi^{-1}J^{(h_1)}:-\tfrac{1}{2}:\Phi^{12}J^{(h_1)}:-\tfrac{1}{2}:\Phi^{12}J^{(h_2)}:-(\tfrac{1}{2}+k)\partial\Phi^{-1}\right.\\
&\left.-\tfrac{1}{2}k\partial\Phi^{12}\right),
\end{align*}
\begin{align*}
L=&\left(-\tfrac{2}{1+2k}\right)\left(:J^{(h_1)}J^{(h_1)}:
-2:J^{(h_1)}J^{(h_2)}:+\left(\tfrac{1+2k}{8}\right):\Phi^{-1}\partial \Phi^{12}:\right. \\
& \left.+(-\frac{1}{2}):\Phi^{12}J^{(e_1)}:+\left(-\tfrac{1+2k}{8}\right):\partial \Phi^{-1}\Phi^{12}:+ \left(\tfrac{1+2k}{16}\right):\partial \Phi^{12}\Phi^{12}:\right.\\
&\left.+(-1)\partial J^{(h_1)}+\left(\tfrac{1+2k}{2}\right)\partial J^{(h_2)}\right),
\end{align*}
\begin{align*}
W=&\tfrac{2\sqrt{1-2k}\sqrt{5+8k}}{5+18k+16k^{2}}\left(
\tfrac{(2+8 k)}{1-2 k}:J^{(h_{1})}J^{(h_{1})}:+\tfrac{(4+16 k)}{1-2 k}:J^{(h_{1})}J^{(h_{2})}:\right.\\
&\left.+\tfrac{(3+6 k)}{1-2 k}:J^{(h_{2})}J^{(h_{2})}:+\tfrac{(3+6 k)}{2-4 k}:\Phi^{-1}\Phi^{12}J^{(h_{1})}:+\tfrac{(3+6 k)}{2-4 k}:\Phi^{-1}\Phi^{12}J^{(h_{2})}:\right.\\
&\left.+\tfrac{(-1-2 k)}{4}:\Phi^{-1}\partial\Phi^{12}:+\tfrac{\left(1+6 k+8 k^2\right)}{4-8 k}:\partial\Phi^{-1}\Phi^{12}:+\tfrac{\left(1+6 k+8 k^2\right)}{-8+16 k}:\partial\Phi^{12}\Phi^{12}:\right.\\
&\left.-\tfrac{2 k (1+4 k)}{-1+2 k}\partial J^{(h_{1})}+\tfrac{\left(1+6 k+8 k^2\right)}{2-4 k}\partial J^{(h_2)}\right),
\end{align*}
\begin{align*}
U=&\left(\tfrac{12 i \sqrt{1+3 k} }{\sqrt{1-2 k} \sqrt{5+8 k} \sqrt{1+5 k+6 k^2}}\right)\left(-:\Phi^{-1}J^{(h_{1})}J^{(h_{1})}:-:\Phi^{-1}J^{(h_{1})}J^{(h_{2})}:\right.\\
&\left.+\tfrac{(1+2 k)}{16}:\Phi^{-1}\partial\Phi^{12}\Phi^{12}:+\tfrac{(1-2 k)}{6}:\Phi^{-1}\partial J^{(h_{1})}:-\tfrac{1}{2}:\Phi^{12}J^{(h_{1})}J^{(h_{2})}:\right.\\
&\left.-\tfrac{1}{2} :\Phi^{12}J^{(h_{2})}J^{(h_{2})}:+\tfrac{(-1-4 k)}{12}:\Phi^{12}\partial J^{(h_{1})}:+\tfrac{(-1-4 k)}{12}:\Phi^{12}\partial J^{(h_{2})}:\right.\\
&\left.+\tfrac{(1+4 k)}{6}:\partial\Phi^{-1}J^{(h_{1})}:+\left(\tfrac{1}{2}+k\right):\partial\Phi^{-1}J^{(h_{2})}:+\tfrac{(1+2 k)}{8}:\partial\Phi^{-1}\Phi^{-1}\Phi^{12}:\right.\\
&\left.+\tfrac{(-1-4 k)}{12}:\partial\Phi^{12}J^{(h_{1})}:+\left(-\tfrac{1}{12}-\tfrac{5 k}{6}\right) :\partial\Phi^{12}J^{(h_{2})}:+\tfrac{\left(1+6 k+8 k^2\right)}{24}\partial^{2}\Phi^{-1}\right.\\
&\left.-\tfrac{k (1+4 k)}{12}\partial^{2}\Phi^{12}\right).
\end{align*}
\emph{Screening operators}
\noindent

Recently was proved how to compute the screening operators of each $W$-algebra obtained via a quantum Hamiltonian reduction. The case considered in the present paper corresponds to the special case when the Cartan subalgebra $\mathfrak{h}$ coincides with $\mathfrak{g}_{0}$ \cite[Theorem 3.9]{Genra2016}. 

Rescaling each current $h\in V_{k+\frac{1}{2}}(\mathfrak{h})$ by $\nu=1/\sqrt{k+1/2}$, we can identified $V_{k+\frac{1}{2}}(\mathfrak{h})$ with the Heisenberg vertex algebra $V_{1}(\mathfrak{h})$ associated to $\mathfrak{h}$. Here we identified $\mathfrak{h}^{\ast}$ with $\mathfrak{h}$ via the bilinear form.
Denote by $M_{\alpha}$ the $V_{1}(\mathfrak{h})$-module with highest weight $\alpha\in \mathfrak{h}^{\ast}$ and highest weight vector $\left|\alpha\right>$;  and denote $\Gamma_{\alpha}(z)$ the intertwiner whose modes are maps $V_{1}(\mathfrak{h})=M_{0}\rightarrow M_{\alpha}$:
\begin{align*}
\Gamma_{\alpha}(z)=s_{\alpha}z^{\alpha_{(0)}}\operatorname{exp}\left(-\sum_{n<0}\tfrac{\alpha_{(n)}}{n}z^{-n}\right)\operatorname{exp}\left(-\sum_{n>0}\tfrac{\alpha_{(n)}}{n}z^{-n}\right),
\end{align*} 
where the map $s_{\alpha}$ is defined by
\begin{align*}
s_{\alpha}\left|0\right>=\left|\alpha\right>, \quad [s_{\alpha},\beta_{(n)}]=0 \quad \forall \; n \neq 0,\; \beta \in \mathfrak{h}^{\ast}.
\end{align*}
Define the screening operators:
\begin{gather*}
Q_{1}=:\Phi_{-1}\Gamma_{\alpha_{1}/\nu}:, \;   Q_{2}=\Gamma_{-\alpha_{2}/\nu}.
\end{gather*}
Then we have
\begin{align}\label{intersectionofscreening}
W_{k}(\mathfrak{g},x,f) \simeq \operatorname{Ker}{Q_{1}}_{(0)} \bigcap \operatorname{Ker}{Q_{2}}_{(0)} \subset V_{1}(h)\otimes F(\mathfrak{g}_{1/2}).
\end{align}
\noindent
Moreover, the image of the isomorphism (\ref{intersectionofscreening}) recovers the free field realization of $W_{k}(\mathfrak{g},x,f)$ inside $V_{k}(\mathfrak{h})\otimes F(\mathfrak{g}_{1/2})$ that we have obtained above.
\bibliographystyle{abbrv}
\bibliography{bibliografia}

\end{document}